\documentclass[12pt]{article}
\usepackage{amssymb,amsfonts,amsmath,amsthm}
\newcommand{\odin}{\mbox {\bf 1}}

\newcommand{\fracd}[2]{\frac {\displaystyle #1}{\displaystyle #2 }}

\newcommand{\nn}{{\mathbb N}}

\newcommand{\rr}{{\mathbb R}}

\newcommand{\calc}{{\mathcal C}}

\newcommand{\calf}{{\mathcal F}}

\newcommand{\veps}{\varepsilon}

\newcommand{\calu}{{\mathcal U}}
\newcommand{\calv}{{\mathcal V}}
\newcommand{\beq}{\begin{eqnarray*}}
\newcommand{\feq}{\end{eqnarray*}}
\newcommand{\beqn}{\begin{eqnarray}}
\newcommand{\feqn}{\end{eqnarray}}

\newtheorem{theorem}{Theorem}
\makeatletter \@addtoreset{equation}{section}\makeatother

\makeatletter \@addtoreset{theorem}{section}\makeatother

\newtheorem{lemma}[theorem]{Lemma}

\newtheorem*{theorem*}{Theorem}
\newtheorem*{remark*}{Remark}
\newtheorem{proposition}[theorem]{Proposition}
\newtheorem{corollary}[theorem]{Corollary}

\title{Optimal control of a stochastic network driven by a fractional Brownian motion input}
\author{Arka P.~Ghosh\thanks{Research supported by National Science Foundation grant DMS-0608634}
\and Alexander Roitershtein
\and Ananda Weerasinghe\thanks{Research supported by  US Army Research Offce grant W911NF0710424.}
\\
$\mbox{}$
\\
\and{Iowa State University}}
\date{August 5, 2008}
\begin{document}
\maketitle
\begin{abstract}
We consider a stochastic control model driven by a fractional
Brownian motion. This model is a formal approximation to a queueing
network with an ON-OFF input process. We study stochastic control
problems associated with the long-run average cost, the infinite
horizon discounted cost, and the finite horizon cost. In addition,
we find a solution to a constrained minimization problem as an
application of our solution to the long-run average cost problem. We
also establish Abelian limit relationships among the value functions
of the above control problems.
\end{abstract}
{\em MSC2000: } primary 60K25, 68M20, 90B22; secondary 90B18.\\
\noindent {\em Keywords: stochastic control, controlled
queueing networks, heavy traffic analysis, fractional Brownian
motion, self-similarity, long-range dependence.}
\section{Introduction}
\label{intro} Self-similarity and long-range dependence of the
underlying data are two important features observed in the
statistical analysis of high-speed communication networks in heavy
traffic, such as local area networks (LAN) (see for instance
\cite{slook, leland, paxson, sahinoglu, taqqu-et-al, will-survey,
will-book} and references therein). In theoretical models such
traffic behavior has been successfully described by stochastic
models associated with fractional Brownian motion (fBM) (see
\cite{hers97, hers98, klin, taqqu-et-al, whitt98, whitt-book}). It
is well known that fBM exhibits both of these features when the
associated Hurst parameter is above $\frac{1}{2}.$ Therefore,
understanding the behavior and control of these stochastic models
are of significant interest. The non-Markovian nature of the
fractional Brownian motion makes it quite difficult to study
stochastic control problems for a state process driven by fBM. The
techniques such as dynamic programming and analysis of the
corresponding Hamilton-Jacobi-Bellmann equations which are the
commonly used tools in the analysis of the stochastic control
problems associated with the ordinary Brownian motion are not
available for fBM-models.
\par
In this paper, we study several basic stochastic control problems for a
queueing model with an input described by a fractional Brownian
motion process. Similar queueing models, but not in the context of control of the state process, were considered for instance
in \cite{conslaw, norros, GZ}. We are aware of only a few solvable stochastic control problems in fBM setting.
Usually, the controlled state process is a solution of a linear (semi-linear in \cite{duncan}) stochastic
differential equation driven by fBM, and the control typically effects the drift term of the SDE.
In particular, the linear-quadratic regulator control problem is addressed
\cite{hu-zhou, esaim} and a stochastic maximum principle is developed and applied to
several stochastic control problems in \cite{MPfbm}. We refer to \cite{hu-zhou}
and to Chapter~9 of the recent book \cite{4fbm-book} for further examples
of control problems in this setting. In contrast to the models considered in the above references, the model described here is
motivated by queueing applications and involves processes with state-constraints. At the end of this section, we discuss
an example of a queueing network which leads to our model. Our analysis relies on a coupling of
the state process with its stationary version (see \cite{conslaw, GZ})
which enables us to address control problems in a non-Markovian setting, and
our techniques are different from those employed in \cite{MPfbm, duncan, hu-zhou, esaim}.
\par
A real-valued stochastic process $W_H=\bigl(W_H(t)\bigr)_{t\geq 0}$
is called fractional Brownian motion with Hurst parameter $H\in
(0,1)$ if $W_H(0)=0$ and $W_H$ is a continuous zero-mean Gaussian
process with stationary increments and covariance function given by
\beq
\mbox{Cov}\bigl(W_H(s),W_H(t)\bigr)=\fracd{1}{2}[t^{2H}+s^{2H}-|t-s|^{2H}],\qquad
s\geq 0, \, t\geq 0. \feq The fractional Brownian motion is a
self-similar process with index $H,$ that is for any $a>0$ the
process $\frac{1}{a^H}\bigl(W_H(at)\bigr)_{t\geq 0}$ has the same
distribution as $\bigl(W_H(t)\bigr)_{t\geq 0}.$ If $H=\frac{1}{2}$
then $W_H$ is an ordinary Brownian motion, and if $H\in
[\frac{1}{2},1)$ then the increments of the process are positively
correlated and the process exhibits long-range dependence, which
means that \beq
        \sum_{n=1}^\infty \mbox{Cov} \bigl(W_H(1),W_H(n+1)-W_H(n)\bigr)=\infty.
\feq
Notice that fBM is a recurrent process, and in particular,
for any $u>0,$ $\lim\limits_{t\to\infty} W_H(t)/t=0$ a.s. and consequently
$\lim\limits_{t\to\infty}\bigl(W_H(t)-ut\bigr)=-\infty$ a.s.
For additional properties and a more detailed description of
the process we refer to \cite{mandvn, nualart, somrt94, shiryaev}.
\par
We consider a single server stochastic processing network having deterministic
service process with rate $\mu>0.$ For any time $t\geq 0,$ the
cumulative work input to the system over the time interval $[0,t]$
is given by $\lambda t+W_H(t),$ where $\lambda$ is a fixed constant
and $W_H$ is a fractional Brownian motion with Hurst parameter $H\in
[\frac{1}{2},1).$ We assume that the service rate $\mu$ satisfies
$\mu>\lambda$ and that the parameter $\mu$ can be controlled. The
workload present in the system at time $t\geq 0$ is given by
$X_x^u(t)$ which is defined in \eqref{cspace} and \eqref{local}
below. Here $x\geq 0$ is the initial workload and $u=\mu-\lambda>0$
is the control variable. Assuming for simplicity that $x=0,$ an
equivalent representation for the process $X_x^u$ is given by (see
\eqref{newx} below for the general case) \beqn \label{process}
X_x^u(t)=\bigl(W_H(t)-ut\bigr)-\inf_{s\in [0,t]}
\bigl(W_H(t)-ut\bigr), \qquad t\geq 0.\feqn For a given arrival
process $W_H,$ this is a common formulation for a simple stochastic
network where the server works continuously unless there are no
customers in the system. The first term above represents the
difference between the cumulative number of job arrivals and
completed services  in the time interval $[0,t]$, and the last term
ensures that the the queue-length is non-negative, and it is a
non-decreasing process which increases only when the queue-length
process is zero. For more examples of such formulations for queueing
networks or stochastic processing networks, we refer to
\cite{harrison-book} and \cite{whitt-book}. The above queueing model
with $W_H$ in \eqref{process} being a fractional Brownian motion was
considered by Zeevi and Glynn in \cite{GZ}. In particular,  we are
motivated by their work.
\par
Our goal is to address several stochastic control
problems related to the control of the workload process $X_x^u$
described above. The organization of the paper is as follows. We
will conclude this section with a motivating example of a queueing
network which leads to our model. In Section~\ref{model}, we
introduce the model and describe three basic stochastic control
problems associated with it, namely the long-run average cost
problem, the infinite horizon discounted cost problem, and the
finite horizon control problem. Here we also discuss some properties
of the reflection map which will be used in our analysis.
\par
In Section~\ref{s-longrun}, we study the long-run average cost
problem. Here we obtain an explicit deterministic representation of
the cost functional for each control $u>0.$ This enables us to
reduce the stochastic control problem to a deterministic
minimization problem. We obtain an optimal control $u^*>0$ and show
its uniqueness under additional convexity assumptions for the
associated cost functions. We show that the value function and the
corresponding optimal control are independent of the initial data.
It is well known that the above property is true for the classical
long-run average cost problem associated with non-degenerate
diffusion processes. Here we show that it remains valid for our
model driven by fBM which is highly non-Markovian. The main results
of this section are given in Theorems~\ref{i-structure} and
\ref{optimalc}. Their generalizations are given in
Theorems~\ref{gen-thm} and \ref{analogue}. Our proofs here rely on
the use of a coupling of the underlying stochastic process with its
stationary version introduced in \cite{conslaw}. In particular, we
show that the coupling time has finite moments in
Proposition~\ref{crutial}. Because of the highly non-Markovian
character of the fractional Brownian motion (it is well known that
fBM cannot be represented as a function of a finite-dimensional
Markov process), coupling arguments in general do not work for the
models associated with fBM (we refer to \cite{hairer} for an
exception). In our case, the coupling is available due to the
uniqueness results related to the reflection map described in
Section~\ref{model}.
\par
We use our results in Section~\ref{s-longrun} to obtain an optimal
strategy for a constrained optimization problem in
Theorem~\ref{conth} of Section~\ref{s-constrained}. Similar
stochastic control problems for systems driven by an ordinary
Brownian motion were previously considered in \cite{ata, weera}. An
interesting application of this model to wireless networks is
discussed in \cite{ata}. Our constrained optimization problem (in
the fBM setting) is a basic example of a general class of
problems with an added bounded variation control process in the
model. This class of problems has important applications to the
control of queueing networks, but in fBM setting, it seems to be an
unexplored area of research.
\par
In Section~\ref{s-horizon}, we address the infinite horizon
discounted cost problem associated with a similar cost structure. An
optimal control for this problem is given in Theorem~\ref{vidix}.
\par
In Section~\ref{abelian}, we establish Abelian limit relationships
among the value functions of the three stochastic control problems
introduced in Section~\ref{model}. The main result of this section
is stated in Theorem~\ref{aben}. We show that the long term
asymptotic for the finite horizon control problem and the asymptotic for the infinite
horizon discounted control problem, as the discount factor
approaches zero, share a common limit. This limit turns out to be
the value of the long-run average cost control problem. Our proof
also shows that the optimal control for the discounted cost problem
converges to the optimal control for the long-run average cost
control problem when the discount factor approaches zero. A similar
result holds also for the optimal control of the finite horizon
problem, as the time horizon tends to infinity. For a class of
controlled diffusion processes analogous results were previously
obtained in \cite{weera}.
\paragraph{Motivating example.} We conclude
this section with a description of a queueing network
related to the internet traffic data in  which the weak limit of a
suitably scaled queue-length process satisfies \eqref{cspace} and
\eqref{local} (which are reduced to the above \eqref{process} when
the initial workload $x$ equals 0). For more details on this model
we refer to \cite{taqqu-et-al} and Sections~7 and 8 of
\cite{whitt-book}.
\par We begin by  defining a sequence of queueing networks with state
dependent arrival and service rates,
indexed by an integer $n \geq 1$ and a non-negative real-valued parameter  $\tau \geq
0$. For each $n\geq 1$ and $\tau\geq 0,$ the $(n, \tau)$-th network has only one
server and one buffer of infinite size, and the arrivals and
departures from the system are given as follows. There are $n$ input
sources (e.g. $n$ users connected to the server), and job
requirements of each user is given by the so-called \emph{ON-OFF
process} $\{X_i^{n, \tau}, i \ge 1\}$ as defined in
\cite{taqqu-et-al}, namely each user stays connected to the server
for a random ON-period of time with distribution function $F_{1}$, and stays off during a random
OFF-period of time with distribution
function $F_{2}$. The distribution $F_i$ is assumed to have finite
mean $m_i$ but infinite variance, and in particular, \beq 1-F_i(x)
\sim c_i x^{-\alpha_i}, \feq where $1<\alpha_i<2$ for $i=1,2$. While
connected to the server, each user demands service at unit rate
(sends data-packets at a unit rate to the server for processing).
The server is processing users requests at a constant rate, say
$\mu_{n,\tau}$. Assume that the ON and OFF-periods are all
independent (for each user as well as across users), the ON-OFF
processes have  stationary increments, and average rate of arrival
of jobs (packets) from each source (customer) is given by $\lambda =
m_1/(m_1+m_2)$. The queue-length at time $t\ge 0$ is given by
\[X^{n,\tau}(t) = X^{n, \tau}(0) + \sum_{i=1}^n \int_0^t X_i^{n,
\tau}(s) ds - \mu_{n,\tau} t + L^{n, \tau}(t), \] where $L^{n,
\tau}$ is a non-decreasing process that starts from 0, increases
only when $X^{n, \tau}$ is zero, and ensures that $X^{n, \tau}$ is
always non-negative. Physically, this implies that the server is
non-idling, i.e it serves jobs continuously as long as the buffer is
non-empty. The second term in the right-hand side of the above equation
represents cumulative number of packets sent to the server by all
the $n$ customers in the interval $[0,t]$. We will assume that $X^{n,
\tau}(0)=x_{n, \tau},$ where $x_{n, \tau}$ are fixed non-negative real numbers
for each $n$ and $\tau.$ In this setup, $\tau >0$ represents the time scaling parameter,
and it is well known (see \cite{funda} or Theorem~7.2.5 of \cite{whitt-book}) that
\[\tau^{-H} n^{-\frac{1}{2}}\sum_{i=1}^n \int_0^{\tau \cdot} X_i^{n,
\tau}(s) ds \; \; \Rightarrow \; \;  W_H(\cdot),\] when $n
\rightarrow \infty$ first and then $\tau \rightarrow \infty$. Here
$W_H$ is a fractional Brownian motion with Hurst parameter  $H=
(3-\min\{\alpha_1,\alpha_2\})/2\in \bigl(\frac{1}{2},1\bigr)$, and the convergence is the weak
convergence in the space $D\bigl([0,\infty),\rr\bigr)$ of
right-continuous real functions with left limits equipped with the
standard $J_{\alpha,1}$ topology (see \cite{whitt-book} for details).
\par
From the above convergence result it can be derived (see
Theorem~8.7.1 in \cite{whitt-book}) that if the service rates $\mu_{n,\tau}$ satisfy the
following heavy traffic assumption\[\tau^{-H}
n^{-\frac{1}{2}}\sum_{i=1}^n \left(\mu_{n\tau}-n\lambda \right)\tau
\, \to \, u, \] as $(n,\tau) \to \infty$, then a
suitably scaled queue-length satisfies equations \eqref{cspace} and
\eqref{local}. More precisely, if the above heavy traffic condition
is satisfied, and $\tau^{-H} n^{-\frac{1}{2}} x_{n,\tau} \rightarrow
x$, then the scaled queue-length $\tau^{-H} n^{-\frac{1}{2}}
X^{n,\tau}(\tau \cdot)$ converges weakly to a limiting process
$X^u_x(\cdot)$ that satisfies \eqref{cspace}, \eqref{local} if we let
$n \rightarrow \infty$ first and then $\tau \rightarrow \infty.$
\par
Hence, we see that with a ``super-imposed" ON-OFF input source and
deterministic services times for the queueing processes, a
suitably-scaled queue length in the limit satisfies our model. With
a cost structure similar to either \eqref{ergodic} or
\eqref{reduced-j}, or \eqref{ifunct} for the queueing network
problem,  one can consider the problem described in this paper as a
formal fractional Brownian control problem (fBCP) of the
corresponding control problem for the queueing network. We, however,
do not attempt to solve the queueing control problem in this paper.
A solution the limiting control problem provides useful insights
into the queueing network control problem (see for instance
\cite{harri3}). For a broad class of queueing problems, it has been
shown that the value function of the Brownian control problem (BCP)
is a lower bound for the minimum cost in the queueing network
control problem (see \cite{ghosh-budhi-2}). In many situations, the
solution to the BCP can be utilized to obtain optimal strategies for
the queueing network control problem (cf. \cite{bellwill},
\cite{ghosh-budhi-1}, \cite{GhWe2} etc.). Here, we study just the
Brownian control problem, which is an important problem in its own
right. Our explicit solution to the fBCP can be considered as an
``approximate solution" to the queueing network problem.
\section{Basic setup}
\label{model} In this section we define the controlled state process
(Section~\ref{ss-mintro}), describe three standard control problems
associated with it (Section~\ref{copr}), and also discuss some basic
properties of a reflection mapping which is involved in the
definition of the state process (Section~\ref{ss-skorokhod}).
\subsection{Model}
\label{ss-mintro} Let $\bigl(W_H(t)\bigr)_{t\geq 0}$ be a fractional
Brownian motion with Hurst parameter $H\geq 1/2$ and let
$\sigma(\cdot)$ be a deterministic continuous function defined on
$[0,\infty)$ and taking positive values. For a given initial value
$x\geq 0$ and a control variable $u\geq 0,$ the controlled state
process $X_x^u$ is defined by \beqn \label{cspace}
X_x^u(t)=x-ut+\sigma(u)W_H(t)+L_x^u(t), \qquad t\geq 0,\feqn where
the process $L_x^u$ is given by \beqn \label{local}
L_x^u(t)=-\min\bigl\{0,\min_{s\in[0,t]}\bigl(x-us+\sigma(u)W_H(s)\bigr)\bigr\},\qquad
t\geq 0. \feqn The control variable $u\geq 0$ remains
fixed throughout the evolution of the state process $X_x^u.$ It
follows from \eqref{cspace} and \eqref{local} that $X_x^u(t)\geq 0$
for all $t\geq 0.$ Notice that the process $L_x^u$ has continuous paths, and it
increases at times when $X_x^u(t)=0.$
\par
The process $X_x^u$ represents the workload process of a single
server controlled queue fed by a fractional Brownian motion, as
described in the previous section (see also \cite{GZ}). For a chosen
control $u\geq 0$ that remains fixed for all $t\geq 0,$ the controller
is faced with a cost structure consisting of the following three additive
components during a time interval $[t,t+dt]:$
\begin{enumerate}
\item a control cost $h(u)dt,$
\item a state dependent holding cost $C\bigl(X_x^u\bigr)dt,$
\item a penalty of $p\,dL_x^u(t),$ if the workload in the system is empty.
\end{enumerate}
Here $p\geq 0$ is a constant, and $h$ and $C$ are non-negative
continuous functions satisfying the following basic assumptions:
\begin{itemize}
\item [(i)] The function $h$ is defined on $[0,+\infty),$ and
\beqn \label{h-function} \mbox{$h$ is non-decreasing and continuous,
$h(0)\geq 0,$ $\lim\limits_{u\to +\infty}h(u)=+\infty.$} \feqn
\item [(ii)] The functions  $C$ is also defined on $[0,+\infty),$ and it is a non-negative, non-decreasing
continuous function which satisfies the following polynomial growth
condition: \beqn \label{growth} 0\leq C(x)\leq K(1+x^\gamma) \feqn
for some  positive constants $K>0$ and $\gamma>0.$
\end{itemize}
We will sometimes assume convexity of $h$ and $C$ in order
to obtain sharper results, such as the uniqueness of the optimal
controls.
\subsection{Three control problems}
\label{copr} Here we formulate three cost minimization problems for
our model. In the {\em long-run average cost minimization problem}
(it is also called {\em ergodic control problem}), the controller
minimizes the cost functional \beqn \nonumber
I(u,x)&:=&\limsup\limits_{T\to\infty}\fracd{1}{T}E\Bigl(\int_0^T
\bigl[h(u)+C\bigl(X_x^u(t)\bigr)\bigr]dt+\int_0^T p\,dL_x^u(t)\Bigr)\\
\label{ergodic} &=& h(u)+\limsup\limits_{T
\to\infty}\fracd{1}{T}E\Bigl(\int_0^T
C\bigl(X_x^u(t)\bigr)dt+p\,L_x^u(T)\Bigr), \feqn subject to the
constraint $u>0$ for a fixed initial value $x\geq 0.$ Notice that since
$L_x^u(t)$ is an increasing process, the integral with respect to $L_x^u(t)$
can be defined an ordinary Riemann-Stieltjes integral. The value function of this problem
is given by \beqn \label{value} V_0(x)=\inf_{u>0}I(u,x). \feqn
\par
In Section~\ref{s-longrun}, we show that $I(u,x)$ and hence also the
value function $V_0(x)$ are actually independent of $x.$ In
addition, we show the existence of a finite optimal control $u^*>0$
and also prove that $u^*$ is unique if the functions $h$ and $C$ are
convex. We apply the results on the long run average cost problem to
find an optimal strategy for a constrained optimization problem, in
Section~\ref{s-constrained}.
\par
In Section~\ref{s-horizon}, we solve the infinite horizon
discounted cost minimization problem for the case when
$\sigma(u)\equiv 1$ in \eqref{cspace}. In this problem it is assumed
that the controller wants to minimize the cost functional \beqn
\label{reduced-j} J_\alpha(u,x):=E\Bigl(\int_0^\infty e^{-\alpha t}
\bigl[ h(u)+C\bigl(X_x^u(t)\bigr)\bigr]dt+p\int_0^\infty e^{-\alpha
t}dL_x^u(t)\Bigr). \feqn subject to $u>0$ for a fixed initial value
$x\geq 0.$ Here the discount factor $\alpha>0$ is a strictly
positive constant. The value function in this case is given by \beqn
\label{valued} V_\alpha(x)=\inf_{u>0}J_\alpha(u,x). \feqn
We study the asymptotic behavior of this model in Section~\ref{abelian}.
When $\alpha$ approaches zero, we prove that $\lim\limits_{\alpha\to
0^+} \alpha J_\alpha(u,x)=I(u,x)$ for any control $u>0.$ Furthermore, we show that
$\lim\limits_{\alpha\to 0^+} \alpha V_\alpha(x)=V_0(x)$ and the
optimal controls for the discounted cost problem converges to that
of the long-run average cost problem as $\alpha$ tends to zero.
In Section~\ref{abelian}, we also consider the {\em finite horizon control problem} with
the value function $V(x,T)$ defined by \beqn\label{fvalue}
V(x,T):=\inf_{u>0} I(u,x,T), \feqn where \beqn \nonumber
I(u,x,T)&:=&E\Bigl(\int_0^T
\bigl[h(u)+C\bigl(X_x^u(t)\bigr)\bigr]dt+ pL_x^u(T)\Bigr)\\
\label{ifunct} &=&h(u)T+p\,E\bigl(L_x^u(T)\bigr)+E\Bigl(\int_0^T
C\bigl(X_x^u(t)\bigr)dt\Bigr), \feqn and $p\geq 0$ is a non-negative
constant. We prove that $\lim\limits_{T\to \infty}
\frac{V(x,T)}{T}=V_0(x).$ Furthermore, we show that the optimal
controls for the finite horizon problem converges to that of the
long-run average cost problem, as $T$ tends to infinity.
\subsection{The reflection map}
\label{ss-skorokhod} The model equations \eqref{cspace} and
\eqref{local} have an equivalent representation which is given below
in \eqref{newx} by using the reflection map. Therefore we
briefly discuss some basic properties of the reflection map and of
the representation \eqref{newx}.
\par
Let $\calc([0,\infty),\rr)$ be the space of continuous functions with
domain $[0,\infty).$ The standard reflection mapping
$\Gamma:\calc([0,\infty),\rr)\to \calc([0,\infty),\rr)$ is defined by \beqn
\label{mapping}
\Gamma\bigl(f\bigr)(t)=f(t)+\sup_{s\in[0,t]}(-f(s))^+, \feqn for
$f\in \calc([0,\infty),\rr).$ Here and henceforth we use the notation
$a^+:=\max\{0,a\}.$ This mapping is also known as the {\em Skorokhod map}
or the {\em regulator map} in different contexts. For a
detailed discussion we refer to \cite{kruk, whitt-book}.
\par
In our model \eqref{cspace}-\eqref{local}, we can write $X_x^u$ as
follows: \beqn \label{newx} X_x^u(t)=\Gamma\bigl(x-u
e+\sigma(u)W_H\bigr)(t), \feqn where $e(t):=t$ for $t\geq 0,$ is the
identity map.
\par
Note that according to the definition, $\Gamma\bigl(f\bigr)(t)\geq
0$ for $t\geq 0.$ We will also use the following two standard facts
about $\Gamma$ (see for instance \cite{kruk, whitt-book}). First,
we have \beqn \label{ineq} \sup_{t\in[0,T]}
|\Gamma\bigl(f\bigr)(t)|\leq 2\sup_{t\in[0,T]} |f(t)|, \feqn for
$f\in \calc([0,\infty),\rr).$ Secondly, let $f$ and $g$ be two
functions in $\calc([0,\infty),\rr)$ such that $f(0)=g(0)$ and
$h(t):=f(t)-g(t)$ is a non-negative non-decreasing function in
$\calc([0,\infty),\rr).$ Then \beqn \label{monot}
\Gamma\bigl(f\bigr)(t)\geq \Gamma\bigl(g\bigr)(t),\qquad \mbox{for
all}~t\geq 0. \feqn We shall also rely on the following convexity
property of the reflected mapping. Let $\alpha\in(0,1)$ and $f$ and
$g$ be two functions in $\calc([0,\infty),\rr).$ Then, $\alpha
f+(1-\alpha)g\in \calc([0,\infty),\rr),$ and \beqn \label{convex}
\Gamma\bigl(\alpha f+(1-\alpha)g\bigr)(t)\leq \alpha
\Gamma\bigl(f\bigr)(t)+(1-\alpha)\Gamma\bigl(g\bigr)(t) \feqn for
all $t\geq 0.$ The proof is straightforward. Let
$F(x)=x^-:=\max\{0,-x\}.$ Then $F$ is a convex function and
therefore $(\alpha f(s)+(1-\alpha)g(s))^-\leq \alpha
f^-(s)+(1-\alpha)g^-(s).$ Consequently, $\sup_{s\in [0,t]} (\alpha
f(s)+(1-\alpha)g(s))^- \leq \alpha \sup_{s\in
[0,t]}f^-(s)+(1-\alpha)\sup_{s\in [0,t]}g^-(s).$ Since $\sup_{s\in
[0,t]}(-f(s))^+=\sup_{s\in [0,t]}f^-(s),$ the inequality
\eqref{convex} follows from the definition \eqref{mapping}.\par
The reflection map also satisfies the following minimality property: 
if $\psi, \eta \in \calc\bigl([0,\infty),\rr\bigr)$ are such that 
$\psi$ is non-negative, $\eta(0)=0,$ $\eta$ is non-decreasing, 
and $\psi(t)=\varphi(t)+\eta(t)$ for $t\geq 0,$ then 
\beqn \label{maxip}
\psi(t)\geq \Gamma(\varphi)(t)\quad \mbox{and}\quad  
\eta(t)\geq \sup_{s\in [0,t]} \bigl(-\varphi(s)\bigr)^+,\qquad \mbox{for all}~t \geq 0.
\feqn
\section{Long-run average cost minimization problem}
\label{s-longrun} In this section we address the control problem
defined in \eqref{ergodic}-\eqref{value}. First we find a solution
to the control problem for the particular case when $\sigma(u)\equiv
1$ in \eqref{cspace}. This is accomplished in Sections~\ref{ss-simple}--\ref{ss-exists}.
In Section~\ref{extent}, we show that the general case can be
reduced to this simplified version.
\subsection{Reduction of the cost structure}
\label{ss-simple} 
The controlled state space process $X_x^u$ (corresponding to $\sigma(u)\equiv 1$)
has the form \beqn \label{simple} X_x^u(t)=x-ut+W_H(t)+L_x^u(t), \qquad
t\geq 0. \feqn The following lemma simplifies the expression for the
cost functional \eqref{ergodic} by computing
$\lim\limits_{T\to\infty} \frac{1}{T} E\bigl(L_x^u(T)\bigr).$
\par
\begin{lemma}
\label{help} Let $X_x^u$ be given by \eqref{simple}. Then
$\lim\limits_{T\to\infty} \frac{1}{T} E\bigl(L_x^u(T)\bigr)=u.$
\end{lemma}
\begin{proof}
Since $u>0,$ using \eqref{newx}, \eqref{ineq}, and \eqref{monot}, we
obtain: \beq 0\leq X_x^u(t)\leq
X_x^0(t)=\Gamma\bigl(x+W_H\bigr)(t)\leq
2\bigl(|x|+\sup_{s\in[0,t]}|W_H(s)|\bigr). \feq By the
self-similarity of fractional Brownian motion process, \beq E\bigl(\sup_{s\in[0,T]}
|W_H(s)|\bigr)\leq K_1T^H, \feq where $K_1\in (0,\infty)$ is a constant
independent of $T$ (see for instance \cite[p. ~296]{nualart}). Therefore, \beqn \label{proof3.1} 0\leq
E\bigl(X_x^u(T)\bigr)\leq 2K_1\bigl(|x|+T^H\bigr). \feqn
Consequently, $\lim\limits_{T\to\infty} \frac{1}{T}
E\bigl(X_x^u(T)\bigr)=0.$ Since $W_H(T)$ is a mean-zero Gaussian
process, it follows from \eqref{simple} that $\frac{1}{T}
E\bigl(L_x^u(T)\bigr)-u=\frac{1}{T}
\Bigl(E\bigl(X_x^u(T)\bigr)-x\Bigr).$ Letting $T$ tend to infinity
completes the proof of the lemma.
\end{proof}
\begin{remark*}
\label{afterproof} Lemma~\ref{help} with literally the same proof as
above remains valid if $X_x^u$ satisfies \eqref{cspace} instead of
\eqref{simple}.
\end{remark*}
With the above lemma in hand, we can represent the cost functional
\eqref{ergodic} and reformulate the long-run average cost
minimization problem as follows. The controller minimizes \beqn
\label{newmodel}
I(u,x)=\bigl(h(u)+pu\bigr)+\limsup\limits_{T\to\infty} \fracd{1}{T}
E\Bigl(\int_0^T C\bigl(X_x^u(t)\bigr)dt\Bigr) \feqn subject to $u>0$
and $X_x^u$ given in \eqref{simple}. Note that the above reduction shows that the original minimization
problem \eqref{ergodic} reduces to the case $p=0$ with the function
$h(u)$ replaced by $h(u)+pu.$
\par
Our next step is to analyze the cost component
$\limsup\limits_{T\to\infty} \fracd{1}{T} E\Bigl(\int_0^t
C\bigl(X_x^u(t)\bigr)dt\Bigr).$ The following results are described in \cite{conslaw} and \cite{overflow},
and are collected in \cite{GZ} in a convenient form for our application. We summarize
them here using our notation.
\begin{itemize}
\item [(i)] The random sequence $X_0^u(t)$ with $t\geq 0$ converges weakly, as $t$ goes to infinity, to
the random variable \beqn \label{zuidof} Z_u:=\max_{s\geq 0}\{W_H(s)-us\}. \feqn
\item [(ii)] There is a probability space supporting the processes $X_0^u,$ $L_0^u$ (and hence
$X_x^u$ as well as $L_x^u$ for any $x\geq 0$) and a stationary
process $X_u^*=\{X_u^*(t):t\geq 0\}$ such that \beqn
\label{pathwise} X_u^*(t)=W_H(t)-ut+\max\{X_u^*(0),L_0^u(t)\},\qquad
t\geq 0, \feqn and \beqn \label{in-distr} X_u^*(t)\overset{D}{=}Z_u, \qquad t\geq 0, \feqn where
$``{\overset{D}{=}}"$ denotes the equality in distribution and $Z_u$ is defined in \eqref{zuidof}.
\item [(iii)]  The tail of the stationary distribution satisfies \beqn
\label{gz-estimate} \lim_{{\mathfrak z}\to\infty} {\mathfrak z}^{2H-2}\log P(Z_u\geq {\mathfrak z})=-\theta^*(u),  \feqn
where $\theta^*(u)$ is given by \beqn \label{theta}
\theta^*(u)=\fracd{u^{2H}}{2H^{2H}(1-H)^{2(1-H)}}>0.
\feqn
In particular, all the moments of $Z_u$ are finite.
\end{itemize}
Throughout the rest of the paper, we use this probability space where
all these processes are defined. Using \eqref{simple} and \eqref{local}, we can write for $t\geq 0,$
\beqn \label{lagrange} X_x^u(t)=W_H(t)-ut+\max\{x,L_0^u(t)\}, \feqn
where \beqn \label{localtime} L_0^u(t)=-\inf_{s\in
[0,t]}\bigl(W_H(s)-us\bigr)=\sup_{s\in [0,t]}\bigl(us-W_H(s)\bigr).
\feqn The above representation \eqref{lagrange} and
\eqref{localtime} agrees with \eqref{pathwise} if the process
$X_x^u$ is initialized with $X_u^*(0).$
\subsection{A coupling time}
\label{coupling}
The following coupling argument is crucial to address the optimal
control problems. In particular, it enables us to deal with the last term 
of $I(u,x)$ in \eqref{newmodel}.
\begin{proposition}
\label{crutial} Let $u>0$ and the initial point $x\geq 0$ be fixed.
Consider the state process $X_x^u$ in \eqref{simple} and the stationary process
$X_u^*$ of \eqref{pathwise} and \eqref{in-distr}. Then the following
results hold:
\begin{itemize}
\item [(i)] There is a finite stopping time $\tau_0$ such that $X_x^u(t)=X_u^*(t)$
for all $t\geq \tau_0.$ Furthermore, $E(\tau_0^\beta)<\infty$ for
all $\beta \geq 0.$
\item [(ii)] The cost functional $I(u,x)$ defined in \eqref{ergodic} is finite and independent of $x,$
that is $I(u,x)=I(u,0)<\infty$ for $x\geq 0.$ Consequently, the
value function $V_0(x)=\inf\limits_{u>0} I(u,x)$ is also finite and
independent of $x,$ that is $V_0(x)=V_0(0)<\infty$ for $x\geq 0.$
\end{itemize}
\end{proposition}
\begin{proof}
For $y>0$ introduce the stopping time \beq
\lambda_y=\inf\{t>0:L_x^u(t)>x+y\}. \feq The stopping time
$\lambda_y$ is finite a.s. because $\lim\limits_{t\to+\infty}
L_x^u(t)\geq \lim\limits_{t\to+\infty}\bigl(ut-W_H(t)\bigr)=+\infty$ a.s.
Define the stopping time $\tau_0$ by \beqn
\label{taust} \tau_0=\inf\{t>0:L_0^u(t)>x+X_u^*(0)\}. \feqn Here
$X_u^*$ is the stationary process which satisfies \eqref{pathwise}
and \eqref{in-distr}. That is $\tau_0=\lambda_{X_u^*(0)}$ a.s.
It follows that for $t\geq \tau_0,$ we have $L_0^u(t)\geq
L_0^u(\tau_0)=x+X_u^*(0)$ and $X_u^*(0)\geq 0.$ Therefore, it
follows from \eqref{lagrange} and \eqref{localtime} that
$X_x^u(t)=W_H(t)-ut+L_u^0(t)=X_u^*(t)$ for $t\geq \tau_0.$
\par
Next, we prove that $E(\tau_0^\beta)<+\infty$ for each $\beta\geq
0.$ Clearly, without loss of generality we can assume that
$\beta\geq 1.$ We then have: \beqn \label{tau-estimate}
E(\tau_0^\beta) \leq \sum_{m=0}^\infty E\bigl(\lambda_{m+1}^\beta
\cdot \odin_{[m\leq X^*_u(0)<m+1]}\bigr) \leq \sum_{m=0}^\infty
\bigl[E (\lambda_{m+1}^{2\beta}) P( X^*_u(0)\geq m) \bigr]^{1/2},
\feqn where in the last step we have used the Cauchy-Schwartz
inequality. Since $X_u^*(0)$ has the same distribution as
$Z_u=\sup_{s\geq 0} \{W_H(s)-us\},$ it follows from
\eqref{gz-estimate} that for all $m$ large enough, \beqn
\label{ngz-estimate} P(X_u^*(0)\geq m)\leq
e^{-\frac{1}{2}\theta^*(u)m^{2(1-H)}},  \feqn  where $\theta^*(u)$
is defined in \eqref{theta}. Next, we estimate
$E(\lambda_m^{2\beta})$ for $m\geq 0$ and $\beta\geq 1.$ For
$m\in\nn,$ let $b_m=x+m$ and $T_m=\frac{2b_m}{u}.$ We have \beqn
\nonumber E(\lambda^{2\beta}_m)&=&2\beta\int_0^\infty
t^{2\beta-1} P(\lambda_m>t)dt\,=\,2\beta\int_0^\infty t^{2\beta-1} P(L_0^u(t)\leq x+m)dt \\
\label{festi}&\leq& T_m^{2\beta}+2\beta\int_{T_m}^\infty t^{2\beta-1}
P(L_0^u(t)\leq b_m)dt, \feqn where the second equality is due to the
fact that $P(\lambda_m>t)=P(L_0^u(t)\leq x+m)$ according to the
definition of $\lambda_m.$ Notice that \beqn \label{sesti}
P\bigl(L_0^u(t)\leq b_m\bigr)&=&P\bigl(\sup_{s\in
[0,t]}\{us-W_H(s)\}\leq b_m\bigr)\,\leq\, P\bigl( W_H(t)\geq
ut-b_m\bigr), \feqn and recall that $Z:=\frac{W_H(t)}{t^H}$ has a
standard normal distribution. Therefore, by \eqref{sesti}, for
$t>T_m$ we have \beqn \nonumber P\bigl(L_0^u(t) \leq b_n\bigr)&\leq&
P\bigl( W_H(t)\geq ut-b_m\bigr) \,\leq\, P\Bigl( W_H(t)\geq \frac{ut}{2}\Bigr)\\
\label{testi} &=& P\Bigl(Z\geq \frac{ut^{1-H}}{2}\Bigr). \feqn It follows from
\eqref{festi} and \eqref{testi} that \beqn \nonumber
E(\lambda^{2\beta}_m)&\leq &
T_m^{2\beta}+2\beta\int_0^\infty t^{2\beta -1} P\Bigl(Z\geq \frac{ut^{1-H}}{2}\Bigr)dt\\
\nonumber &=& T_m^{2\beta}+2\beta\int_0^\infty t^{2\beta-1}
P\Bigl[\Bigl(\fracd{2Z}{u}\Bigr)^{\frac{1}{1-H}}\geq t\Bigr]dt\\
&=& \label{vesti}
T_m^{2\beta}+E\Bigl[\Bigl(\fracd{2|Z|}{u}\Bigr)^{\frac{2\beta}{1-H}}\Bigr]=
\fracd{4^\beta}{u^{2\beta}}(x+m)^{2\beta}+E\Bigl[\Bigl(\fracd{2|Z|}{u}\Bigr)^{\frac{2\beta}{1-H}}\Bigr]<\infty.
\feqn  The estimates \eqref{ngz-estimate} and \eqref{vesti} imply
that the infinite series in the right-hand side of
\eqref{tau-estimate} converges. Thus $E(\tau_0^\beta)<\infty$ for
all $\beta\geq 1,$ and hence for all $\beta \geq 0.$ This completes
the proof of the first part of the proposition.
\par
We turn now to the proof of part (ii). First, we will prove that
\beqn \label{frefer} E\Bigl(\int_0^\infty
\bigl|C\bigl(X_x^u(t)\bigr)-C\bigl(X_u^*(t)\bigr)\bigr|dt\Bigr)<\infty.
\feqn We will show later that part (ii) of the proposition is a
rather direct consequence of this inequality. \par Notice that \beq
E\Bigl(\int_0^\infty
\bigl|C\bigl(X_x^u(t)\bigr)-C\bigl(X_u^*(t)\bigr)\bigr|dt\Bigl)=E\Bigl(\int_0^{\tau_0}
\bigl|C\bigl(X_x^u(t)\bigr)-C\bigl(X_u^*(t)\bigr)\bigr|dt\Bigl),
\feq where $\tau_0$ is given in \eqref{taust} and $X_u^*$ is given
in \eqref{pathwise} and \eqref{in-distr}. The definition of $\tau_0$
implies $L_0^u(\tau_0)\leq x+X_u^*(0).$ Therefore, it follows from
\eqref{zuidof} and \eqref{pathwise} that for $t\in [0,\tau_0],$ \beq
\max \{X_x^u(t),X_u^*(t)\}\leq Z_u+x+X^*_u(0). \feq Since $C$ is a
non-decreasing function, this implies \beq \max
\bigl\{C\bigl(X_x^u(t)\bigr),C\bigl(X_u^*(t)\bigr)\}\leq
C\bigl(Z_u+x+X^*_u(0)\bigr),\feq and consequently, using the
Cauchy-Schwartz inequality, \beq E\Bigl(\int_0^{\tau_0}
\bigl|C\bigl(X_x^u(t)\bigr)-C\bigl(X_u^*(t)\bigr)\bigr|dt\Bigl)
&\leq& E\Bigl(\tau_0 C\bigl(Z_u+x+X^*_u(0)\bigr)\Bigr) \\
&\leq& \Bigl[ E(\tau_0^2) E\Bigl(
\bigl[C\bigl(Z_u+x+X^*_u(0)\bigr)\bigr]^2\Bigr)\Bigr]^{1/2}.\feq
Since $E(\tau_0^2)<\infty$ by part (i) of the lemma, \eqref{frefer}
will follow once we show that \beqn \label{refer1} E\Bigl(
\bigl[C\bigl(Z_u+x+X^*_u(0)\bigr)\bigr]^2\Bigr)<\infty. \feqn Recall
that $X_u^*(0)$ and $Z_u$ have the same distribution. The tail asymptotic
\eqref{gz-estimate} implies that any moment of $Z_u$ is finite.
This fact together with \eqref{growth} yield \eqref{refer1}.
\par
We will now deduce part (ii) of the proposition from \eqref{frefer}.
Toward this end, first observe that, since $X_u^*(t)$ is a stationary process,
\beq \fracd{1}{T} E\Bigl(\int_0^T
C\bigl(X_u^*(t)\bigr)dt\Bigr)= \fracd{1}{T} E\Bigl(\int_0^T
C\bigl(X_u^*(0)\bigr)dt\Bigr)=E\bigl(C(Z_u)\bigr),\feq and recall
that $E\bigl(C(Z_u)\bigr)<\infty$ by \eqref{gz-estimate}.
Then notice that by \eqref{frefer}, \beq && \limsup_{T\to\infty}
\Bigl| \fracd{1}{T} E\Bigl(\int_0^T
C\bigl(X_u^*(t)\bigr)dt\Bigr)-\fracd{1}{T} E\Bigl(\int_0^T
C\bigl(X_x^u(t)\bigr)dt\Bigr)\Bigr| \\ && \qquad \qquad \leq
\limsup_{T\to\infty} \fracd{1}{T} E\Bigl(\int_0^\infty
\bigl|C\bigl(X_x^u(t)\bigr)-C\bigl(X_u^*(t)\bigr)\bigr|dt\Bigr)=0.\feq
Therefore \beqn \label{part} \lim\limits_{T\to\infty}\frac{1}{T} E\bigl(\int_0^T
C\bigl(X_x^u(t)\bigr)dt\bigr)= E\bigl(C(Z_u)\bigr),\feqn which implies
part (ii) of the proposition in view of \eqref{newmodel}. Therefore, the proof
of the proposition is complete.
\end{proof}
\begin{remark*}
The above proposition is in agreement with a result in Theorem~1 of
\cite{GZ} which shows that the asymptotic distributions, as $t$
approaches infinity, for $M(t)= \max_{s\in [0,t]} X_x^u(s)$ and
$M^*(t)= \max_{s\in [0,t]} X_u^*(s)$ coincide.
\end{remark*}
\subsection{Properties of $E\bigl(C(Z_u)\bigr).$}
\label{ss-expectation} For $u>0,$ let $G(u)=E\bigl(C(Z_u)\bigr).$
We are interested in the behavior of $G(u)$ in view of the identity \eqref{part}.
\begin{lemma}
\label{g-expo} Let $G(u)=E\bigl(C(Z_u)\bigr)$ where $Z_u$ is defined
in \eqref{zuidof}. Then the following results hold:
\begin{itemize}
\item [(i)] $G(u)$ is a decreasing and continuous function of $u$ on $[0,\infty).$
\item [(ii)] If $C(x)$ is a convex function then so is $G(u).$
\item [(iii)] $\lim\limits_{u\to 0^+} G(u)=+\infty.$
\end{itemize}
\end{lemma}
\begin{proof}
First we observe that the polynomial bound \eqref{growth} on the
growth of $C$ combined with \eqref{gz-estimate}, which
describes the tail behavior of $Z_u,$ imply that
$G(u)=E\bigl(C(Z_u)\bigr)$ is finite for each $u\geq 0.$ It is a
decreasing function of $u$ because $C$ is non-decreasing while
$Z_{u_1}\leq Z_{u_2}$ if $u_1>u_2.$
\par
To complete the proof of part (i), it remains to show that $G(u)$ is a continuous function. To
prove this, first notice that, according to the definition \eqref{zuidof},
$Z_u$ is a continuous function of the variable $u$ a.s., as shown below. Indeed, if $t_u$ is a random time
such that $Z_u=W_H(t_u)-ut_u$ and $u\in (0,v),$ then
\beq
Z_u \geq Z_v \geq W_H(t_u)- vt_u=Z_u- t_u(v-u).
\feq
Hence $\lim_{v\to u^+} Z_v=Z_u.$ A similar argument shows that $\lim_{v\to u^-} Z_v=Z_u.$
Therefore, since $C$ is a continuous function, $C(Z_u)$ is a continuous function of $u$ with
probability one. Since $C(Z_u)$ is monotone in $u,$ the dominated convergence theorem implies 
the continuity of $G.$ 
\par To prove part (ii), fix constants $r\in [0,1],$ $u_1>0,$ $u_2>0,$ and let $\bar
u_r=ru_1+(1-r)u_2.$ Then $Z_{\bar u_r}=\sup_{t\geq 0} \{W_H(t)-\bar
u_r t\}\leq rZ_{u_1}+(1-r)Z_{u_2}.$ If $C$ is a non-decreasing
convex function, we have $E\bigl(C(Z_{\bar u_r})\bigr)\leq
rE\bigl(C(Z_{u_1})\bigr)+(1-r)E\bigl(C(Z_{u_2})\bigr).$ Hence $G$ is
convex, and the proof of part (ii) is complete.
\par
Turning to the part (iii), we first notice that $Z_0=\sup_{t\geq 0}
W_H(t)=+\infty$ with probability one. Let $(u_n)_{n\geq 0}$ be any
sequence monotonically decreasing to zero. Then, $Z_{u_n}$ is increasing and
hence there exists a limit (finite or infinite)
$\lim\limits_{n\to\infty} Z_{u_n}=L$ and $Z_{u_n}\leq L$
for all $n\geq 0.$ Thus $W_H(t)-u_nt\leq L$ for all $n\geq 0$ and
$t\geq 0.$ By letting $n$ go to infinity we obtain $\sup_{t\geq 0}
W_H(t)\leq L,$ and consequently $L=+\infty$ with probability one.
Therefore $\lim\limits_{u\to 0^+} Z(u)=+\infty$ a.s. Since $C$ is a
non-decreasing function, the monotone convergence theorem implies
that $\lim\limits_{u\to 0^+} E\bigl(C(Z_u)\bigr)=+\infty.$ This
completes the proof of the lemma.
\end{proof}
\subsection{Existence of an optimal control}
\label{ss-exists} In the following two theorems we provide a
representation of the cost functional $I(u,x)$ as well as the
existence and uniqueness results for the optimal control $u^*>0.$
\begin{theorem}
\label{i-structure}
Let $I(u,x)$ be the cost functional of the long-run average cost problem described in \eqref{newmodel}. Then
\item [(i)] $I(u,x)$ is independent of $x$ and has the representation
\beqn
\label{i-repr}
I(u):=I(u,x)=
h(u)+pu+G(u),
\feqn
where $G(u)$ is given in Lemma~\ref{g-expo}.
Furthermore, $I(u)$ is finite for each $u>0$ and is continuous in $u>0.$
\item [(ii)] $\lim\limits_{u\to 0^+} I(u)=+\infty$ and $\lim\limits_{u\to\infty} I(u)=+\infty.$
\item [(iii)] If $h(x)$ and $C(x)$ are convex functions, then $I(u)$ is also convex.
\end{theorem}
\begin{proof}
Part (i) follows from \eqref{newmodel}, Proposition~\ref{crutial}, and Lemma~\ref{g-expo}.
\par
The first part of claim (ii) follows from the fact that $I(u)\geq
G(u)$ along with part~(iii) of Lemma~\ref{g-expo}. To verify the
second part, notice that $I(u)\geq h(u)$ for all $u>0,$ and
$\lim\limits_{u\to +\infty} h(u)=+\infty.$ Consequently,
$\lim\limits_{u\to+\infty} I(u)=+\infty.$
\par
Part (iii) follows from the representation \eqref{i-repr} combined with the part~(ii) of Lemma~\ref{g-expo}. This completes
the proof of the theorem.
\end{proof}
\begin{theorem}
\label{optimalc}
\item [(i)] There is an optimal control $u^*>0$ such that for all $x\geq 0$ we have
\beq I(u^*)=\min_{u>0} I(u,x) \feq where $I$ is given in
\eqref{i-repr}. In particular, $u^*$ is independent of $x.$
\item [(ii)] In the case $p>0,$ if $h$ and $C$ are convex functions, then $u^*$ is unique.
\item [(iii)] In the case $p=0,$ if $h$ is a strictly convex function and $C$ is a convex function, then $u^*$ is unique.
\end{theorem}
\begin{proof}
Since $I(u)$ is a continuous function, part (i) follows from parts
(i) and (ii) of Theorem~\ref{i-structure}.
\par
If $h$ and $C$ are convex functions, the representation \eqref{i-repr} yields that $I$
is a strictly convex function when $p>0.$ Therefore, $u^*$ is unique in this case.
\par
In the case $p=0,$ if $h$ is a strictly convex function and $C$ is a convex function, the results
follows from the representation \eqref{i-repr} in a similar way as in the case $p>0.$
\end{proof}
\subsection{Generalizations}
\label{extent} In this section we generalize the results in
Theorem~\ref{i-structure} and \ref{optimalc} to the more general
model introduced in \eqref{cspace}. Notice that for a given control
$u>0$ fixed in  \eqref{cspace}, the self-similarity of fBM yields
that the process $\widehat W_H$ defined by \beq \widehat
W_H(t)=\sigma(u)W_H\Bigl(\frac{t}{\sigma(u)^{\frac{1}{H}}}\Bigr),~t\in\rr,
\feq is a fractional Brownian motion. Let \beq
Y_x^u(t)=X_x^u\Bigl(\frac{t}{\sigma(u)^{\frac{1}{H}}}\Bigr).\feq
Then $Y_x^u$ satisfies \beqn \label{dual}
Y_x^u(t)=x-\frac{ut}{\sigma(u)^{\frac{1}{H}}}+\widehat
W_H(t)+\widehat L_x^u(t), \feqn where \beq \widehat
L_x^u(t)=L_x^u\Bigl(\frac{t}{\sigma(u)^{\frac{1}{H}}}\Bigr).\feq
Using \eqref{local} and change of the variable
$s=\frac{t}{\sigma(u)^{\frac{1}{H}}},$ we observe that \beqn
\label{elhat} \widehat
L_x^u(t)=L_x^u\Bigl(\frac{t}{\sigma(u)^{\frac{1}{H}}}\Bigr)=
-\min\Bigl\{0,\min_{s\in[0,t]}\Bigl(x-\frac{su}{\sigma(u)^{\frac{1}{H}}}+\widehat
W_H(s)\Bigr)\Bigr\}. \feqn The equations \eqref{dual} and \eqref{elhat} are analogous
to \eqref{cspace} and \eqref{local}.
\par
We next consider the change in the cost structure due to the change
of the variable $s=\frac{t}{\sigma(u)^{\frac{1}{H}}}.$ We notice
that \beq \fracd{1}{T}\int_0^T C\bigl(X_x^u(t)\bigr)dt=
\fracd{1}{M(T)}\int_0^{M(T)} C\bigl(Y_x^u(t)\bigr)dt, \feq where
$M(T)=\sigma(u)^{\frac{1}{H}}T.$ Therefore, using the results in
Theorem~\ref{i-structure}, we obtain \beqn \label{3.29}
\lim\limits_{T\to\infty} \fracd{1}{T}E\Bigl(\int_0^T
C\bigl(X_x^u(t)\bigr)dt\Bigr)=
G\Bigl(\frac{u}{\sigma(u)^{\frac{1}{H}}}\Bigr). \feqn We have the
following result.
\begin{theorem}
\label{gen-thm} Consider the controlled state process $X_x^u$
defined by \eqref{cspace} and \eqref{local} with the cost functional
$I(x,u)$ given in \eqref{ergodic}. Define $f:[0,\infty)\to
[0,\infty)$ by \beqn \label{f-function}
f(u)=\fracd{u}{\sigma(u)^{\frac{1}{H}}}. \feqn Then \beqn
\label{i-identity} I(u,x)=h(u)+p\,u+G\bigl(f(u)\bigr), \feqn where
the function $G$ is given in Lemma~\ref{g-expo}. Furthermore,
$I(u,x)$ is independent of $x$ (we will henceforth denote the cost
function by $I(u)$).
\end{theorem}
\begin{proof}
The same argument as in the proof of Lemma~\ref{help} yields
$\lim\limits_{T\to\infty} \frac{1}{T}E\bigl(L_x^u(T)\bigr)=u.$
Combining this result with \eqref{3.29} we obtain the representation
\eqref{i-identity}.
\end{proof}
Our next result is analogous to Theorems~\ref{i-structure} and \ref{optimalc}.
\begin{theorem}
\label{analogue} Assume that the function $f$ in \eqref{f-function}
is continuous and $\lim\limits_{u\to 0^+} f(u)=0.$ Then, with
$I(u)=I(u,x)$ as in \eqref{i-identity}, we have:
\begin{itemize}
\item [(i)] $\lim\limits_{u\to 0^+} I(u)=\lim\limits_{u\to +\infty} I(u)=+\infty,$ and
$I(u)$ is a finite continuous function on $[0,\infty).$ Furthermore, there is
a constant $u^*>0$ such that $I(u^*)=\min_{u>0} I(u).$
\item [(ii)] If $f$ is a concave increasing function then the statements
similar to parts (i) and (ii) of Theorem~\ref{optimalc} (regarding
the uniqueness $u^*$) hold.
\end{itemize}
\end{theorem}
The proof of this theorem is a straightforward
modification of the proofs of Theorems~\ref{i-structure} and \ref{optimalc} and therefore is omitted.
\begin{remark*}
\label{sec3}
One can further generalize our model to
cover the following situation. For given positive continuous functions $b(u)$ and $\sigma(u)$ let
\beq
X_x^u(t)=x+\sigma(u)W_h(t)-b(u)t+L_x^u(t),
\feq
where for $u>0,$
\beq
L_x^u(t)=-\min\left\{0,\min_{s\in[0,t]}\bigl(x-b(u)s+\sigma(u)W_H(s)\bigr)\right\}.
\feq
The optimization problem here is to minimize the cost functional
$I(u,x)$ defined in \eqref{ergodic}.
\par
Following the time change method described in Section~\ref{extent},
one can obtain an analogue of Theorem~\ref{analogue} regarding the
derivation of the optimal control. In this situation, the function
$f$ defined in \eqref{f-function} needs to be replaced by
$f(u)=b(u)(\sigma(u))^{-\frac{1}{H}}$ with the assumptions that $f$ is continuous and
$\lim\limits_{u\to0^+} f(u)=0.$ We omit the details of the proof.
\end{remark*}
\section{A constrained minimization problem}
\label{s-constrained} In this section, we address a constrained
minimization problem that can be solved by using our results in
Section~\ref{s-longrun}. Our model here is of the form \beqn
\label{modelhere} Y^u_x(t)=x-ut+ \sigma(u)W_H(t)+K^u(t), \feqn where
$\sigma$ is a non-negative continuous function, $K^u(\cdot)$ is a
non-negative non-decreasing right-continuous with left limits (RCLL)
process adapted to the natural filtration $(\calf_t)_{t\geq 0},$
where $\calf_t$ is the $\sigma$-algebra generated by $\{W_H(s):0\leq
s\leq t\}$ augmented with all the null sets. Furthermore, $K^u(0)=0$
and the process $K^u$ is chosen by the controller in such a way that
the state process $Y^u$ is constrained to non-negative reals. In
this situation, the controller is equipped with two controls: the
choice of $u>0$ and the choice of $K^u$ process subject
non-negativity of the $Y^u$ process.
\par
Throughout this section we keep the initial state $x\geq 0$ fixed.
We will deduce using the results of the previous section that the value of this
minimization problem as well as the optimal control are not affected
by the initial data.
\par
Let $m>0$ be any fixed positive constant. The constrained
minimization problem we would like to address here is the
following: \beqn
\nonumber
\mbox{Minimize}\quad&&
\\
\label{aim} && \limsup\limits_{T\to\infty} \fracd{1}{T} E\Bigl(
\int_0^T \bigr[h(u)+C\bigl(Y^u(t)\bigr)\bigr] dt\Bigr)
\\
\nonumber
\mbox{Subject to:}&&
\\
\label{constraint} && \limsup\limits_{T\to\infty}
\fracd{E\bigl(K^u(T)\bigr)}{T}\leq m. \feqn A controlled
optimization problem of this nature for diffusion processes was
considered in \cite{ata}, and for a more complete treatment in the
case of diffusion processes we refer to \cite{weera}.
\par
Fix any integer $m>0$ and define a class of state processes
$\calu_m$ as follows: \beq \calu_m =\Bigl\{(Y^u,K^u):~Y^u(t)\geq
0~\mbox{for}~t\geq 0, ~\mbox{\eqref{modelhere} is satisfied},
~\limsup\limits_{T\to\infty}\fracd{E\bigl(K^u(T)\bigr)}{T}\leq
m\Bigr\}. \feq From our results in Section~\ref{s-longrun} it
follows that for any $u\leq m,$ the pair $(X_x^u,L_x^u)$ in
\eqref{cspace} and \eqref{local} belongs to $\calu_m,$ and hence
$\calu_m$ is non-empty. Therefore, the constrained minimization
problem is to find \beq \inf_{(Y^u,K^u)\in \calu_m}
\limsup\limits_{T\to\infty} \fracd{1}{T} E\Bigl(\int_0^T
\bigr[h(u)+C\bigl(Y^u(t)\bigr)\bigr]dt\Bigr). \feq In this section
we make the following additional assumptions:
\begin{itemize}
\item [(i)] For functions $h$ and $C$ we assume:
\beqn
\mbox{$h$ is strictly convex and satisfies
\eqref{h-function},\, $C$ is convex and satisfies \eqref{growth}.}
\label{extra1}
\feqn
\item [(ii)] Let $f(u)=\frac{u}{\sigma(u)^{\frac{1}{H}}}.$ Then
\beqn \mbox{$ f(u)>0,$~$\lim\limits_{u\to 0^+} f(u)=0,$ and $f$ is a
convex increasing function.} \label{extra2} \feqn
\end{itemize}
The following lemma enables us to reduce the set $\calu_m$
to the collection of processes $(X_x^u,L_x^u)$ described in the previous section, with
$u\leq m.$
\begin{lemma}
Let $u>0$ and let $(Y^u,K^u)$ be a pair of processes satisfying \eqref{modelhere}.
Consider $(X_x^u,L_x^u)$ which satisfies \eqref{cspace}, \eqref{local}, and is defined 
on the same probability space as $(Y^u,K^u).$ Then
\item [(i)] $L_x^u(t)\leq K^u(t)$ and $X_x^u(t)\leq Y^u(t)$ for $t\geq 0.$
\item [(ii)] $u=\lim\limits_{T\to\infty} \frac{1}{T} E\bigl(L_x^u(T)\bigr)\leq
\limsup\limits_{T\to\infty}\frac{1}{T} E\bigl(K^u(T)\bigr).$
\item [(iii)] $G\bigl(f(u)\bigr)=\lim\limits_{T\to\infty} \frac{1}{T}
E\bigl(\int_0^T C\bigl(X_x^u(t)\bigr)dt\bigr) \leq \limsup\limits_{T\to\infty} \frac{1}{T}
E\bigl(\int_0^T C\bigl(Y^u(t)\bigr)dt\bigr).$
\end{lemma}
\begin{proof}
Since, $Y^u\geq 0,$ $K^u(0)=0,$ and $K^u$ is non-increasing process, the minimality property of the reflection map 
stated in \eqref{maxip} implies that $L_x^u(t) \leq  K^u(t)$ and $X_x^u(t)\leq Y_x^u(t)$ for $t\geq 0.$
\par
Next, observe that part (ii) of the lemma follows from the result in
part (i) while the identity $u=\lim\limits_{T\to\infty} \frac{1}{T}
E\bigl(L_x^u(T)\bigr)$ is implied by Lemma~\ref{help} (see also
the remark right after the proof of Lemma~\ref{help}).
\par
Finally, part (iii) of the lemma follows from \eqref{3.29}, part (i),
and from the fact that $C$ is non-decreasing. The proof of the lemma is complete.
\end{proof}
Let \beq \calv_m=\{(X_x^u,L_x^u):~\mbox{\eqref{cspace},\eqref{local}
are satisfied and in addition $u\leq m$}\}. \feq From the above
lemma it is clear that \beqn \nonumber && \inf_{(K^u,Y^u)\in\calu_m}
\lim\limits_{T\to\infty} \fracd{1}{T} E\Bigl(\int_0^T
\bigr[h(u)+C\bigl(X_x^u(t)\bigr)\bigr]dt\Bigr)\\  && \qquad \qquad
\qquad = \inf_{(X_x^u,L_x^u)\in \calv_m} \limsup\limits_{T\to\infty}
\frac{1}{T} E\Bigl(\int_0^T
\bigr[h(u)+C\bigl(Y^u(t)\bigr)\bigr]dt\Bigr). \label{dominance}
\feqn Therefore, our minimization problem is reduced. Next, we can
use the results in Section~\ref{extent} and write for any $u>0,$
\beqn \label{use-results} \lim\limits_{T\to\infty} \fracd{1}{T}
E\Bigl(\int_0^T
\bigr[h(u)+C\bigl(X_x^u(t)\bigr)\bigr]dt\Bigr)=h(u)+G\bigl(f(u)\bigr).
\feqn Here $G$ is given by Lemma~\ref{g-expo}, and $f$ is described
in \eqref{f-function} and \eqref{extra2}. Consequently, \beqn
\nonumber && \inf_{(X_x^u,L_x^u)\in \calv_m}
\lim\limits_{T\to\infty} \fracd{1}{T} E\Bigl(\int_0^T
\bigr[h(u)+C\bigl(X_x^u(t)\bigr)\bigr]dt\Bigr) \\  && \qquad \qquad
\qquad = \inf\{h(u)+G\bigl(f(u)\bigr):0<u\leq m\}.
\label{inside-vim} \feqn We next consider the optimal control
described in Theorem~\ref{analogue} corresponding to the case $p=0.$
In virtue of the assumptions \eqref{extra1} and \eqref{extra2}, the
optimal control is unique and we will label it by $u_0^*>0.$ We have
the following result.
\begin{theorem}
\label{conth}
Let
\beqn
\label{u-trancated}
u^*(m)=\left\{
\begin{array}{rcl}
m&\mbox{\em if}&m<u_0^*,\\
u_0^*&\mbox{\em if}&m\geq u_0^*
\end{array}
\right.,
\feqn
where $u_0^*$ is the unique optimal control
in Theorem~\ref{analogue} corresponding to $p=0.$
\par
Then the pair $(X_x^{u^*(m)},L_x^{u^*(m)})$ is an optimal process
for the constrained minimization problem \eqref{aim} and
\eqref{constraint}. Furthermore, the optimal control $u^*(m)$ is a
continuous increasing function of the parameter $m.$
\end{theorem}
\begin{proof}
Let $\Lambda(u)=h(u)+G\bigl(f(u)\bigr)$ where $f$ and $G$ are as in
\eqref{use-results}. Then, by the assumptions \eqref{extra1},
\eqref{extra2} and Theorem~\ref{analogue}, $\Lambda$ is a strictly
convex function which is finite everywhere on $(0,\infty).$
Furthermore, $\lim\limits_{u\to0^+} \Lambda(u)=+\infty$ and
$\lim\limits_{u\to +\infty} \Lambda(u)=+\infty,$ and hence $\Lambda$
has a unique minimum at $u_0^*.$ Therefore, $\Lambda$ is strictly
increasing on $(u_0^*,\infty).$ Clearly, with $u^*(m)$ defined in \eqref{u-trancated},
\beq
\Lambda\bigl(u^*(m)\bigr)=\inf_{u\leq m} \Lambda(u),\feq
and $u^*(m)$ is the unique number which has this property.
By \eqref{dominance} and \eqref{inside-vim}, we have \beq
\Lambda\bigl(u^*(m)\bigr)=\inf_{(X_x^u,L_x^u)\in \calv_m}
\limsup\limits_{T\to\infty} \fracd{1}{T} E\Bigl(\int_0^T
\bigr[h(u)+C\bigl(Y^u(t)\bigr)\bigr]dt\Bigr). \feq Consider the pair
of processes $(X_x^{u^*(m)},L_x^{u^*(m)})$ defined in \eqref{cspace} and
\eqref{local}. Then, in virtue of Lemma~\ref{help}, we have
$\lim\limits_{T\to\infty} \frac{1}{T}
E\bigl(L^{u^*(m)}(T)\bigr)=u^*(m)\leq m,$ and by Theorem~3.6, \beq
\limsup\limits_{T\to\infty} \fracd{1}{T} E\Bigl(\int_0^T
\bigr[h(u)+C\bigl(X_x^{u^*(m)}(t)\bigr)\bigr]dt\Bigr)=
\Lambda\bigl(u^*(m)\bigr). \feq Hence $(X_x^{u^*(m)},L_x^{u^*(m)})$
describes an optimal strategy. This completes the proof of the theorem.
\end{proof}
\begin{remark*}
Notice that the above optimal control $u^*(m)$
is independent of the initial point $x.$
\end{remark*}
\section{Infinite horizon discounted cost minimization problem}
\label{s-horizon} In this section we define an optimal control $u^*$
for the infinite horizon discounted cost functional given in
\eqref{reduced-j}. Throughout this section
we assume that $\sigma(u)\equiv 1,$ the state process $X^u_x$
satisfies \eqref{simple}, and that the functionals $h$ and $C$ are
convex in addition to the assumptions stated in \eqref{h-function}
and \eqref{growth}. In contrast with Section~\ref{s-longrun}, 
our methods here do not readily extend to the case where the function 
$\sigma(u)$ is non constant. 
\par
The discounted cost functional $J_\alpha(x,u)$ is given by \beqn \nonumber
J_\alpha(x,u)&=&E\Bigl(\int_0^\infty e^{-\alpha t}
\bigr[h(u)+C\bigr(X_x^u(t)\bigr)\bigr]dt +\int_0^\infty e^{-\alpha t}p\,dL_x^u(t)
\Bigr)
\\
&=& \label{j-sigma1} \fracd{h(u)}{\alpha}+ E\Bigl(\int_0^\infty
e^{-\alpha t} \bigr[C\bigr(X_x^u(t)\bigr) +\alpha
p\,L_x^u(t)\bigr]dt\Bigr). \feqn Here $\alpha>0$ is a constant
discount function. To derive the last equality above we have used Fubini's
theorem to obtain $\int_0^\infty e^{-\alpha t} dL_x^u(t)=\alpha
\int_0^\infty e^{-\alpha t}L_x^u(t)dt.$
\par
Let \beq J_{\alpha,1}(x,u)=E\Bigl(\int_0^\infty e^{-\alpha
t}C\bigr(X_x^u(t)\bigr)dt\Bigr) \feq and \beqn \label{jey2}
J_{\alpha,2}(x,u)=E\Bigl(\int_0^\infty e^{-\alpha
t}L_x^u(t)dt\Bigr). \feqn Next we use the convexity of the
reflection mapping described in \eqref{convex} to establish the
convexity of the cost functional with respect to $u.$
\begin{lemma}
\label{convexity}
Let $x\geq 0$ be fixed and let $C$ be a convex function satisfying assumption \eqref{growth}. Then,
$J_{\alpha,1}(x,u)$ and $J_{\alpha,2}(x,u)$ introduced above are finite for each $u\geq 0$ and are convex
in $u$-variable.
\end{lemma}
\begin{proof}
By \eqref{proof3.1}, we have the bound $E\bigl(L_x^u(t)\bigr)\leq ut +K_0(1+t^H),$ where $K_0>0$ is a
constant independent of $t.$ This implies, by using Fubini's
theorem, that $J_{\alpha,2}(x,u)$ is finite.
\par
Next, by using \eqref{frefer} we obtain \beq
\Bigl|J_{\alpha,1}(x,u)- E\Bigl(\int_0^\infty e^{-\alpha
t}C\bigr(X_u^*(t)\bigr)dt\Bigr)\Bigr| \leq E\Bigl(\int_0^\infty
\bigl|C\bigl(X_x^u(t)\bigr)-C\bigl(X_u^*(t)\bigr)\bigr|dt\Bigr)<\infty. \feq 
But, using the stationary of $X_u^*,$ we have $E\Bigl(\int_0^\infty e^{-\alpha
t}C\bigr(X_u^*(t)\bigr)dt\Bigr)=\alpha^{-1} E\bigl(C(Z_u)\bigr),$
where $Z_u$ is given in \eqref{zuidof}. Notice that $E\bigl(C(Z_u)\bigr)$ is
finite because $C$ has polynomial growth and in virtue of \eqref{gz-estimate}.
Consequently, $J_{\alpha,1}(x,u)$ is also finite.
\par
To establish convexity of $J_{\alpha,1}(x,u),$ first recall that
$X_x^u=\Gamma\bigl(x+W_H-ue\bigr)$, where $e(t)\equiv t$ for $t\geq
0,$ and $\Gamma$ is the reflection mapping described in
Section~\ref{ss-skorokhod}. Now let $u_1\geq 0,$ $u_2\geq 0,$ and
$r\in (0,1).$ Then, for $t \geq 0,$ \beq &&
x+W_H(t)-\bigl(ru_1+(1-r)u_2\bigr)e(t)\\ && \qquad \qquad \qquad
=r\bigl(x+W_H(t)-u_1e(t)\bigr)+(1-r)\bigl(x+W_H(t)-u_2e(t)\bigr).
\feq Since the reflection map $\Gamma$ satisfies the convexity
property \eqref{convex}, we have: \beqn \label{conx} X_x^{\bar u_r}
(t)\leq rX_x^{u_1} (t)+(1-r)X_x^{u_2} (t), \feqn where $\bar
u_r=ru_1+(1-r)u_2.$ Next, since $C$ is a convex non-decreasing
function, \eqref{conx} implies that $C\bigl(X_x^{\bar u_r}
(t)\bigr)\leq rC\bigl(X_x^{u_1} (t)\bigr)+(1-r)C\bigl(X_x^{u_2}
(t)\bigr)$ for $t\geq 0.$ From this it follows that \beq
J_{\alpha,1}(x,\bar u_r) \leq r
J_{\alpha,1}(x,u_1)+(1-r)J_{\alpha,1}(x,u_2). \feq Hence
$J_{\alpha,1}(x,u)$ is convex in $u$-variable.
\par
Next, by \eqref{local} and \eqref{mapping}, we have \beq
\Gamma\bigl(x+W_H-ue\bigr)(t)-(x+W_H-ue)(t)=L_x^u(t),\qquad t\geq 0.
\feq Then, since $\Gamma$ is convex in $u$-variable by
\eqref{convex}, the process $L_x^u$ is also convex in $u$-variable.
Consequently, with $\bar u_r=r u_1+(1-r)u_2,$ we obtain \beq
L_x^{\bar u_r} (t)\leq rL_x^{u_1} (t)+(1-r)L_x^{u_2} (t), \qquad
t\geq 0.\feq Finally, it is evident that $J_{\alpha,2}(x,u)$ is
convex in $u$-variable from the definition \eqref{jey2}.
\end{proof}
\begin{corollary}
\label{c-convexity}
Under the conditions of Lemma~\ref{convexity}, the discounted cost functional $J_\alpha(x,u)$ is finite for each $u\geq 0$ and is convex
in $u$-variable.
\end{corollary}
\begin{proof}
Notice that
\beqn
\label{jey-repr}
J_\alpha(x,u)=\fracd{h(u)}{\alpha} + J_{\alpha,1}(x,u)+\alpha p J_{\alpha,2}(x,u).
\feqn
By our assumptions, $h$ is a convex function and $p\geq 0$ and $\alpha>0$ are
constants. Therefore, the claim follows from Lemma~\ref{convexity}.
\end{proof}
The above lemma and the corollary lead to the following result.
\begin{theorem}
\label{vidix} Consider the $X_x^u$ process satisfying \eqref{simple}
and the associated discounted cost functional $J_\alpha(x,u)$ described in
\eqref{j-sigma1}. Then, for each initial point $x\geq 0,$ there is
an optimal control $u^*\geq 0$ such that \beq J_\alpha(x,u^*)=\inf_{u\geq
0} J_\alpha(x,u)\equiv V_\alpha(x), \feq where $V_\alpha(x)$ is the value
function of the discounted cost problem defined in \eqref{valued}.
\end{theorem}
\begin{proof}
Fix any $x\geq 0$ and $\alpha>0.$ By Corollary~\ref{c-convexity},
$J_\alpha(x,u)$ is finite for each $u\geq 0$ and is convex in
$u$-variable. By \eqref{jey-repr}, we have that $J_\alpha(x,u)\geq
\frac{h(u)}{\alpha}$ and hence, since $\lim\limits_{u\to \infty}
h(u)=+\infty,$ we obtain $\lim\limits_{u\to\infty}
J_\alpha(x,u)=+\infty.$ Since $J_\alpha(x,u)$ is convex in $u$-variable, we can conclude that there is a $u^*\geq 0$
(which may depend on $x$) such that $J_\alpha(x,u^*)=\inf_{u\geq 0} J_\alpha(x,u).$
This completes the proof of the theorem.
\end{proof}
\begin{remark*}
Notice that in contrast with the long-run average cost minimization
problem, we cannot rule out the possibility $u^*=0$ here.
\end{remark*}
\begin{corollary}
For the special case $p=0,$ assume further that $h(x)$ is constant on an
interval $[0,\delta]$ for some $\delta>0.$ Then, for every initial point
$x\geq 0,$ the optimal control $u^*$ is strictly positive.
\end{corollary}
\begin{proof}
It follows from \eqref{jey-repr} that $J_\alpha(x,u)=J_{\alpha,1}(x,u).$ The function
$C$ is increasing and $X_x^{u_1}(t)< X_x^{u_2}(t)$ for
all $t\geq 0$ and $u_1>u_2.$ Therefore $J_{\alpha,1}(x,u_1)\leq J_{\alpha,1}(x,u_2)$ for
$u_1>u_2.$ Consequently, $J_\alpha(x,0)\geq J_\alpha(x,u)$ for all $u>0.$ Hence,
we can find an optimal control $u^*>0,$ and the proof of the
corollary is complete.
\end{proof}
\section{Finite horizon problem and Abelian limits}
\label{abelian} In this section, we establish Abelian limit
relationships among the value functions of three stochastic control
problems introduced in Section~\ref{copr}. The main result is stated
in Theorem~\ref{aben} below. Throughout this section, for
simplicity, we assume that $\sigma(u)\equiv 1$ in the model described in
\eqref{cspace}.
\par
We begin with the existence of an optimal control for the finite
horizon control problem introduced in Section~\ref{copr}.
\begin{proposition}
For any fixed $x\geq 0$ and $T>0,$ we have:
\item [(i)] $I(u,x,T)$ is finite for each $u>0$ and is continuous in $u>0.$
\item [(ii)] $\lim\limits_{u\to\infty} I(u,x,T)=+\infty.$
\item [(iii)] If $h$ and $C$ are convex functions, then $I(u,x,T)$ is a convex function of
the variable $u.$
\end{proposition}
\begin{corollary} For any fixed $x\geq 0$ and $T>0,$ we have:
\item [(i)] There is an optimal control $u^*(x,T)\geq 0$ such that $I\bigl(u^*(x,T)\bigr)=\min_{u>0} I(u,x,T).$
\item [(ii)] In the case $p>0,$ if $h$ and $C$ are convex functions, then $u^*(x,T)$ is unique.
\item [(iii)] In the case $p=0,$ if $h$ is a strictly convex function and $C$ is a convex function, then $u^*(x,T)$ is unique.
\end{corollary}
The proof of the proposition and of the corollary is a
straightforward adaptation of the corresponding proofs given in
Section~\ref{s-longrun}, and therefore is omitted.
\par
The following theorem is the main result of this section.
\begin{theorem}
\label{aben} 
Let $X_x^u$ satisfy \eqref{simple} and let $V_0,$ $V_\alpha(x),$
$V(x,T)$ be the value functions defined in \eqref{value}, \eqref{valued}, and \eqref{fvalue}
respectively. Then the following Abelian limit relationships hold:
\beq
\lim\limits_{\alpha\to 0^+}\alpha V_\alpha(x)=\lim\limits_{T\to\infty} \fracd{V(x,T)}{T}=V_0. \feq
\end{theorem}
We prove this result in Propositions~\ref{abet} and \ref{abef} below. The following technical lemma
gathers necessary tools to establish $\lim\limits_{\alpha\to 0^+}\alpha V_\alpha(x)=V_0.$
\begin{lemma}
\label{gathel}
Let $u>0$ be given and $X_x^u$ satisfy \eqref{simple}. Consider the cost functional $J_\alpha(x,u)$
as defined in \eqref{reduced-j}. Then
\beq
\lim\limits_{\alpha\to 0^+} \alpha J_\alpha(x,u)=I(u),
\feq
where $I(u)$ is described in \eqref{i-repr}.
\end{lemma}
\begin{proof}
First consider $\lim\limits_{\alpha\to 0^+}\alpha E\bigl(\int_0^\infty e^{-\alpha t} dL_x^u(t)\bigr),$ where
$L_x^u$ is as in \eqref{local}. Similarly to \eqref{j-sigma1}, we have
\beqn
\label{asquare}
\alpha E\Bigl(\int_0^\infty e^{-\alpha t} dL_x^u(t)\Bigr)
&=&\alpha^2 E\Bigl(\int_0^\infty e^{-\alpha t} L_x^u(t)dt\Bigr)
=\alpha^2 \int_0^\infty e^{-\alpha t} E\bigl(L_x^u(t)\bigr)dt,
\feqn
where we used Fubini's theorem to obtain the last identity. By \eqref{simple},
\beq
E\bigl(L_x^u(t)\bigr)=ut+E\bigl(X_x^u(t)\bigr)-x.\feq
Therefore,
\beqn
\nonumber
\alpha^2 \int_0^\infty e^{-\alpha t} E\bigl(L_x^u(t)\bigr)dt&=&
\alpha^2u\int_0^\infty e^{-\alpha t} tdt+
\alpha^2 \int_0^\infty e^{-\alpha t} E\bigl(X_x^u(t)\bigr)dt-\alpha x
\\
\label{afurther}
&=&
u+\alpha^2 \int_0^\infty e^{-\alpha t} E\bigl(X_x^u(t)\bigr)dt -\alpha x.
\feqn
Further, by \eqref{proof3.1}, $0\leq  E\bigl(X_x^u(t)\bigr) \leq K_0(1+t^H),$
where the constant $K_0>0$ is independent of $t$ and $u.$ Thus
\beqn
\nonumber
0&\leq& \alpha^2 \int_0^\infty e^{-\alpha t} E\bigl(X_x^u(t)\bigr)dt
\leq K_0\alpha^2 \int_0^\infty e^{-\alpha t} (1+t^H)dt \\ \label{gammaf} &\leq&
K_0\bigl(\alpha+\mbox{Gamma}(H)\alpha^{1-H}\bigr),
\feqn
where $\mbox{Gamma}(H):=\int_0^\infty e^{-t}t^Hdt=\alpha^{1+H}\int_0^\infty e^{-\alpha t} t^Hdt$
is the Gamma function evaluated at $H.$
\par
Since $H\in (0,1),$ it follows from \eqref{gammaf} that
$\lim\limits_{\alpha\to 0^+} \alpha^2 \int_0^\infty e^{-\alpha t}
E\bigl(X_x^u(t)\bigr)dt=0.$ Hence, using \eqref{asquare} and
\eqref{afurther}, we obtain \beqn \label{flimit}
\lim\limits_{\alpha\to 0^+}\alpha E\Bigl(\int_0^\infty e^{-\alpha t}
dL_x^u(t)\Bigr)=u. \feqn We next consider $\lim\limits_{\alpha\to
0^+}\alpha E\bigl(\int_0^\infty e^{-\alpha t}
C\bigl(X_x^u(t)\bigr)dt\bigr).$ It follows from \eqref{frefer} that
\beq && \Bigl|E\Bigl(\int_0^\infty e^{-\alpha t}
C\bigl(X_x^u(t)\bigr)dt\Bigr)-E\Bigl(\int_0^\infty e^{-\alpha
t}C\bigl(X_u^*(t)\bigr)\bigr|dt\Bigr) \Bigr|\\ && \qquad \qquad \leq
E\Bigl(\int_0^\infty
\bigl|C\bigl(X_x^u(t)\bigr)-C\bigl(X_u^*(t)\bigr)\bigr|dt\Bigr)<\infty,
\feq where $X_u^*$ is the stationary process described in
\eqref{pathwise} and \eqref{in-distr}. Therefore, \beqn \nonumber &&
\lim\limits_{\alpha\to 0^+}\alpha E\bigl(\int_0^\infty e^{-\alpha t}
C\bigl(X_x^u(t)\bigr)dt\bigr)=\lim\limits_{\alpha\to 0^+}\alpha
E\bigl(\int_0^\infty e^{-\alpha t} C\bigl(X_u^*(t)\bigr)dt\bigr)\\
\label{ccomp} && \qquad \qquad = \lim\limits_{\alpha\to 0^+}\alpha
\int_0^\infty e^{-\alpha t} E\bigl(C(Z_u)\bigr)dt=
E\bigl(C(Z_u)\bigr)=G(u),\feqn where $G(u)$ is defined in Section~\ref{ss-expectation}.
It follows from \eqref{reduced-j},
\eqref{flimit}, and \eqref{ccomp} that \beq
\lim\limits_{\alpha\to 0^+}\alpha
J_\alpha(x,u)=h(u)+pu+G(u)=I(u), \feq where $I(u)$ is
given in \eqref{i-repr}. This completes the proof of the lemma.
\end{proof}
The next proposition contains the proof of the first part of Theorem~\ref{aben}.
\begin{proposition}
\label{abet}
Let $X_x^u$ satisfy \eqref{simple} and $V_\alpha(x)$
be the corresponding value function defined in \eqref{valued}. Then
\beq \lim\limits_{\alpha\to 0^+}\alpha V_\alpha(x)=V_0<\infty, \feq
where $V_0$ is the value of the long-run average cost minimization
problem given in \eqref{value}.
\end{proposition}
\begin{proof}
Fix the initial point $x\geq 0.$ For any $u>0,$ we have $\alpha
V_\alpha (x)\leq \alpha J_\alpha (x,u).$ Hence, by
Lemma~\ref{gathel}, $\limsup\limits_{\alpha\to 0^+}\,\alpha
V_\alpha(x) \leq \lim\limits_{\alpha\to 0^+}\alpha
J_\alpha(x,u)=I(u).$ Therefore, minimizing the right-hand side over
$u>0,$ we obtain \beq \limsup\limits_{\alpha\to 0^+}\,\alpha
V_\alpha(x) \leq \inf_{u>0} I(u)=V_0. \feq It remains to show the
validity of the reverse inequality, namely that \beqn \label{lowerb}
\liminf\limits_{\alpha\to 0^+}\,\alpha V_\alpha(x) \geq \inf_{u>0}
I(u)=V_0. \feqn To this end, consider a decreasing to zero sequence
$\alpha_n>0,$ $n\in \nn,$ such that \beqn \label{towarda}
\liminf\limits_{\alpha\to 0^+}\,\alpha
V_\alpha(x)=\lim\limits_{n\to\infty}\,\alpha_n V_{\alpha_n}(x).
\feqn Fix any $\veps_0>0$ and let $u_n>0,$ $n\in\nn,$ be a sequence
such that $V_{\alpha_n}(x)+\veps_0>J_{\alpha_n}(x,u_n).$ Then \beqn
\label{alphan}
\alpha_nV_{\alpha_n}(x)+\alpha_n\veps_0>\alpha_nJ_{\alpha_n}(x,u_n)\geq
h(u_n).\feqn Letting $n\to\infty$ we obtain \beqn \label{vioe}
\limsup\limits_{n\to \infty} h(u_n)\leq \limsup\limits_{n\to
\infty}\alpha_n V_{\alpha_n}(x)\leq V_0.\feqn Since
$\lim\limits_{x\to\infty} h(x)=+\infty,$ this implies that $u_n$ is
a bounded sequence. That is, there is $M>0$ such that $u_n\in (0,M)$
for all $n\in\nn.$ Therefore, without loss of generality we can
assume that $u_n$ converges as $n\to\infty$ to some $u_\infty\in
[0,M]$ (otherwise, we can consider a convergent subsequence of
$u_n$).
\par
Let $\delta \in (u_\infty,\infty).$ Then, \beq
\alpha_nJ_{\alpha_n}(x,u_n) \geq h(u_n)+\alpha_n^2p \int_0^\infty
e^{-\alpha_n t} E\bigl(L_x^{u_n}(t)\bigr)dt+ \alpha_n \int_0^\infty
e^{-\alpha_n t} E\bigl[C\bigl(X_x^\delta(t)\bigr)\bigr]dt \feq Since
$E\bigl(L_x^{u_n}(t)\bigr)= E\bigl(X_x^{u_n}(t)\bigr)+u_nt-x\geq
u_nt-x,$ we obtain \beq \alpha_n^2\int_0^\infty e^{-\alpha_n t}
E\bigl(L_x^{u_n}(t)\bigr)dt\geq u_n-\alpha_nx, \feq and hence \beq
\alpha_nJ_{\alpha_n}(x,u_n) \geq h(u_n)+p\,u_n-p \alpha_nx+ \alpha_n
\int_0^\infty e^{-\alpha_n t}
E\bigl[C\bigl(X_x^\delta(t)\bigr)\bigr]dt. \feq Therefore, letting
$n$ go to infinity and using \eqref{ccomp}, we obtain \beqn
\label{vse} \liminf\limits_{n \to \infty}
\alpha_nJ_{\alpha_n}(x,u_n) \geq
h(u_\infty)+p\,u_\infty+E\bigl[C\bigl(Z_\delta\bigr)\bigr]. \feqn In
particular, one can conclude that $u_\infty>0,$ because otherwise,
letting $\delta$ tend to $0$ and using Lemma~\ref{g-expo}, we would
obtain that $\liminf_{n\to\infty} \alpha_nJ_{\alpha_n}(x,u_n)=+\infty$
which contradicts the fact that $\limsup\limits_{n \to \infty}
\alpha_nJ_{\alpha_n}(x,u_n) \leq V_0<\infty$ according to
\eqref{alphan} and \eqref{vioe}.
\par
Therefore $u_\infty>0.$ Letting $\delta$ in \eqref{vse} tend to $u_\infty$ and using again
the continuity of $G(u)=E\bigl[C\bigl(Z_u\bigr)\bigr]$ proved in Lemma~\ref{g-expo}, we obtain
\beqn
\label{alpham}
\liminf\limits_{n\to\infty} \alpha_nJ_{\alpha_n}(x,u_n) \geq
h(u_\infty)+p\,u_\infty+E\bigl[C\bigl(Z_{u_\infty}\bigr)\bigr]=I(u_\infty) \geq V_0.
\feqn
The inequalities \eqref{alphan}, \eqref{vioe}, and \eqref{alpham} combined together yield
\beq
\lim\limits_{n\to\infty} \alpha_nJ_{\alpha_n}(x,u_n) =I(u_\infty)=V_0,
\feq
which completes the proof of the proposition in view of \eqref{towarda}. 
Notice that \eqref{alpham} implies that $u_\infty$ is an optimal control
for the long-run average cost control problem.
\end{proof}
The following proposition includes the second part of Theorem~\ref{aben}.
\begin{proposition}
\label{abef}
Under the conditions of Theorem~\ref{aben}, we have
\beq
\lim\limits_{T\to\infty} \fracd{V(x,T)}{T}=V_0. 
\feq
\end{proposition}
\begin{proof}
It follows from \eqref{ergodic} and \eqref{ifunct} that for any
$x\geq 0$ and $u>0$ we have \beq \limsup_{T\to\infty}
\fracd{I(u,x,T)}{T}=I(u), \feq where $I(u)$ is given in
\eqref{i-repr}. Therefore $\limsup\limits_{T\to \infty}\,\fracd
{V(x,T)}{T} \leq \limsup\limits_{T \to \infty }
\fracd{I(u,x,T)}{T}=I(u),$ and, minimizing the right-hand side over
$u>0,$ we obtain \beqn \label{fuperb} \limsup\limits_{T\to\infty}\,
\fracd{V(x,T)}{T} \leq \inf_{u>0} I(u)=V_0. \feqn It remains to show
that \beqn \label{fowerb}
\liminf\limits_{T\to\infty}\,\fracd{V(x,T)}{T} \geq V_0. \feqn The
proof of \eqref{fowerb} given below is quite similar to that of
\eqref{lowerb}. Consider a sequence of positive reals $(T_n)_{n\in
\nn}$ such that $\lim\limits_{n\to\infty} T_n=+\infty$ and \beq
\liminf\limits_{T\to\infty}\,\fracd{
V(x,T)}{T}=\lim\limits_{n\to\infty}\,\fracd{V(x,T_n)}{T_n}. \feq Fix
any $\veps_0>0$ and for each $n\in\nn$ chose $u_n>0$ such that
$V(x,T_n)+\veps_0>I(u_n,x,T_n).$ Then, in view of \eqref{ifunct} and
\eqref{fuperb}, we have \beq \limsup\limits_{n\to\infty} h(u_n)\leq
\limsup\limits_{n\to\infty}\fracd{V(x,T_n)}{T_n}\leq V_0<+\infty.\feq
Since $\lim\limits_{x\to\infty} h(x)=+\infty,$ this implies that
$u_n$ is a bounded sequence. That is, there is $M>0$ such that
$u_n\in (0,M)$ for all $n\in\nn.$ Taking a further subsequence if
necessary, we can assume without loss of generality that $u_n$
converges to some $u_\infty\in [0,M],$  as $n\to\infty.$
Furthermore, $u_\infty>0$ because if $u_\infty=0,$ then by
\eqref{part} we obtain \beq V_0\geq \limsup\limits_{n\to\infty}
\fracd{V(x,T_n)}{T_n}\geq \limsup\limits_{n\to\infty}
\fracd{I(u_n,x,T_n)}{T_n}\geq E\bigl(C(Z_\delta)\bigr)\qquad
\mbox{for any}~\delta>0. \feq This is impossible in view of
part~(iii) of Lemma~\ref{g-expo}. Therefore $u_\infty>0.$ Let
$v_1,v_2$ be any numbers such that $0<v_1<u_\infty <v_2.$ Then, by
\eqref{ifunct} and an argument similar to the derivation of
\eqref{vse}, we have \beq V(x,T_n)+\veps_0>I(u_n,x,T_n) \geq
h(v_1)T_n+E\bigl(L_x^{v_1}(T_n)\bigr)+\int_0^{T_n}E\bigl[C\bigl(X_x^{v_2}(t)\bigr)\bigr]dt,
\feq for $n$ large enough. Using Lemma~\ref{help} along with
\eqref{part} we deduce that \beq \liminf\limits_{n\to\infty}
\fracd{V(x,T_n)}{T_n}\geq h(v_1)+p\,v_1+E\bigl(C(Z_{v_2})\bigr).
\feq Since $v_1$ and $v_2$ are arbitrary numbers satisfying the
above inequality constraints and $h(u)$ and
$G(u)=E\bigl(C(Z_u)\bigr)$ are continuous functions, this implies
\beq V_0\geq \liminf\limits_{n\to\infty} \fracd{V(x,T_n)}{T_n}\geq
h(u_\infty)+p\,u_\infty+E\bigl(C(Z_{u_\infty})\bigr)=I(u_\infty)\geq
V_0. \feq The proof of \eqref{fowerb}, and hence the proof of the
proposition is complete. In fact, the above inequality also shows
that $u_\infty$ is an optimal control for the long-run average cost
control problem.
\end{proof}
Propositions~\ref{abet} and \ref{abef} combined together yield Theorem~\ref{aben}.
\begin{remark*} The proof of Propositions~\ref{abet} and \ref{abef} imply the 
following results:  
\item [1.] Let $(\alpha_n)_{n\geq 0}$ be a sequence of positive numbers 
converging to zero and let $u_n$ be an $\veps$-optimal control for 
$V_{\alpha_n}(x)$ in \eqref{valued}. Then the sequence $(u_n)_{n\geq 0}$ is bounded and any 
limit point of $(u_n)_{n\geq 0}$ is an optimal control for $V_0$ defined in \eqref{value}.     
\item [2.] Let $(T_n)_{n\geq 0}$ be a sequence of positive numbers such that 
$\lim\limits_{n\to\infty} T_n =+\infty.$ If $u_n$ be an $\veps$-optimal control for
$V(x,T_n)$ in \eqref{fvalue}, then the sequence $(u_n)_{n\geq 0}$ is bounded and any
limit point of $(u_n)_{n\geq 0}$ is an optimal control for $V_0$ defined in \eqref{value}.   
\end{remark*}
\section*{Acknowledgments:}
The work of Arka Ghosh was partially supported by the National Science Foundation grant DMS-0608634.
The work of Ananda Weerasinghe was partially supported by the US Army Research Office  grant W911NF0710424.
We are grateful for their support.

$\mbox{}$
\\
$\mbox{}$
\\
\begin{tabular}{l}
Arka P.~Ghosh \\
303 Snedecor Hall\\
Department of Statistics\\
Iowa State University \\
Ames, IA 50011, USA \\
apghosh@iastate.edu
\end{tabular}
\hfill $\mbox{}$ \hfill
\begin{tabular}{l}
Alexander Roitershtein\\
420 Carver Hall \\
Department of Mathematics \\
Iowa State University \\
Ames, IA 50011, USA \\
roiterst@iastate.edu
\end{tabular}
\\
$\mbox{}$
\\
$\mbox{}$
\\
$\mbox{}$
\\
\begin{tabular}{l}
Ananda Weerasinghe \\
414 Carver Hall \\
Department of Mathematics \\
Iowa State University \\
Ames, IA 50011, USA \\
e-mail: ananda@iastate.edu
\end{tabular}
\end{document}